\documentclass[12pt]{iopart}

\usepackage{url}            
\usepackage{graphicx}
\usepackage{amssymb}
\usepackage{amsthm}
\usepackage{latexsym}
\usepackage{verbatim}
\usepackage{marginnote}
\usepackage[pagebackref]{hyperref} 		
\hypersetup{
    colorlinks,
    citecolor=blue,
    filecolor=magenta,
    linkcolor=red,
    urlcolor=cyan
}

\theoremstyle{plain}
\newtheorem{theorem}{Theorem}[section]

\newtheorem{lemma}[theorem]{Lemma}

\theoremstyle{definition}

\newtheorem{assumption}[theorem]{Assumption}

\def\clo#1{\overline{#1}}
\def\text#1{\mbox{#1}}

\newcommand{\ra}{\rangle}
\newcommand{\la}{\langle}

\newcommand{\C}{c}

\usepackage[letterpaper,margin=1.15in]{geometry}

\begin{document}

\title[Recovery of Pressure and Wave Speed for Photoacoustic Imaging]{Recovery of Pressure and Wave Speed for Photoacoustic Imaging under a Condition of Relative Uncertainty}

\author{Sebasti\'{a}n Acosta}

\address{Department of Pediatrics, Baylor College of Medicine, Texas, USA \\Predictive Analytics Laboratory, Texas Children's Hospital, USA}
\ead{sebastian.acosta@bcm.edu}
\vspace{10pt}
\begin{indented}
\item[]\today
\end{indented}

\begin{abstract}
In this paper we study the photoacoustic tomography problem for which we seek to recover both the initial state of the pressure field and the wave speed of the medium from knowledge of a single boundary measurement. The goal is to propose practical assumptions to define a set of initial conditions and wave speeds over which uniqueness for this inverse problem is guaranteed. The main result of the paper is that given two sets of wave speeds and pressure profiles, they cannot produce the same acoustic measurements if the relative difference between the wave speeds is much smaller than the relative difference between the pressure profiles. Implications for iterative joint-reconstruction algorithms are discussed. 
\end{abstract}

%
\vspace{2pc}

%
%
%
%

\section{Introduction} \label{sec:intro}

Photoacoustic tomography (PAT) is an emerging imaging modality that combines two types of physical fields. The domain to be imaged is illuminated with a short laser pulse that gets absorbed by the medium. The absorbed energy triggers an expansive pressure wave whose initial amplitude is proportional to the optical absorption coefficient of the tissues within the domain. The pressure waves propagate to the domain's boundary where they are measured and processed to recover the initial state of the pressure field. The advantage of this multiwave modality is that biological tissues exhibit high contrast in optical absorption whereas the acoustic waves carry high resolution. Thus, high contrast and high resolution can be achieved simultaneously which offers great potential for biomedical imaging
\cite{Beard2011,Cox2012,Wang-Anastasio-2011,Wang-2009,Wang2012,Wang-Wu-2007}.  
 
Most of the reconstruction methods for PAT assume that the acoustic properties of the medium are known. For a homogeneous wave speed, explicit formulas are available \cite{Finch2004,Kuchment2008,Kunyansky2007,Kunyansky2011,Haltmeier2014,Natterer2012}. Corrections of analytical formulas, valid for asymptotically small variations of sound speed, were investigated in \cite{Dean-Ben2014,Jin2006}. For heterogeneous media, iterative methods have been designed \cite{Stefanov2009,Qian2011,Huang2013,Belhachmi2016,Chervova2016,Stefanov2017a,Haltmeier2017b,Acosta-Montalto-2016,Acosta-Montalto-2015,Palacios2016b}. These methods also require the precise knowledge of the acoustic parameters. However, in practice the wave speed is not known exactly. For soft biological tissues, variations of the wave speed can be as great as $10 \%$ \cite{Jin2006}. If not accounted for in the reconstruction methods, these variations cause blurring and misplacement of features in the reconstructed image. 

It has been noted that the photoacoustic measurements carry information not only about the absorbed optical energy, but also about the wave speed of the medium. Based on this observation, the main question is whether the initial state of the pressure field $u_{0}$ and the wave speed $\C$ can be recovered simultaneously from a single photoacoustic measurement. In general terms, this question of uniqueness is still open. This problem, which is linear in $u_{0}$ and nonlinear in $\C$, is very challenging. However, some progress has been made. Under certain practical assumptions, if one of the two components in $(u_{0},\C)$ is known, the other can be recovered from boundary measurements. See \cite{Stefanov2013} and references therein. Liu and Uhlmann \cite{Liu2015} gave sufficient conditions to recover both the initial pressure profile $u_{0}$ and sound speed $\C$. More precisely, given another pair $(\tilde{u}_{0},\tilde{\C})$, then $\C^{-2} u_{0} = \tilde{\C}^{-2} \tilde{u}_{0}$ if either $\varphi = \C^{-2} u_{0} - \tilde{\C}^{-2} \tilde{u}_{0}$ is a harmonic function or $\varphi$ is independent of one variable in $\mathbb{R}^3$. Oksanen and Uhlmann studied how a modeling error in the wave speed affects the accuracy of the reconstruction of the pressure \cite{Oksanen2014}. Stefanov and Uhlmann concluded that the linearized version of this problem is unstable in any scale of Sobolev spaces \cite{Stefanov2013b}. Kirsch and Scherzer \cite{Kirsch2012} proposed an approach to simultaneously identifying the optical absorbing density and speed of sound based on a family of sectional photoacoustic illuminations and corresponding measurements. These specific illuminations must be focused on cross-sectional planes and special detectors should be employed to neglect out-of-plane signals. 

Computational studies have also been carried out. Treeby et al. proposed a method to select a sound speed that maximizes the sharpness of the reconstructed image \cite{Treeby2011}. Matthews et al. proposed a joint reconstruction method based on an optimization framework and a low-dimensional parametrization of the sound speed \cite{Matthews2018}. Matthews and Anastasio offered an approach based on combining PAT measurements with ultrasound tomography measurements to estimate the wave speed concurrently with the pressure field \cite{Matthews2017b}. A numerical investigation was also performed by Huang et al. \cite{Huang2016a}. The common conclusion from these numerical studies is that severe ill-conditioning is observed for the optimization-based methods if no regularization terms are incorporated.

Motivated by the physical scenario encountered in the PAT problem, we propose some assumptions to ensure the unique recovery of the initial pressure state and wave speed. These assumptions are stated in Section \ref{Sec:Assumptions}. We show that under those conditions, for each pair $(u_{0},\C)$, there is small \emph{region} of the spaces in which they reside, from which no other pair $(\tilde{u}_{0},\tilde{\C})$ can induce the same acoustic measurements. Unfortunately, this region is not a neighborhood of $(u_{0},\C)$. However, the region does coincide with the conditions implicitly assumed in practical/computational scenarios, namely, that the wave speed $\C$ can be estimated a--priori by a known profile $\tilde{\C}$ with much more relative accuracy than the initial state of the pressure field. This result is stated in precise terms at the end of Section \ref{Sec:Assumptions} and the proof is provided in Section \ref{Sec:Proof}. Some final remarks concerning iterative reconstruction algorithms are offered in Section \ref{Sec:Final}.

\section{Assumptions and Main Result} \label{Sec:Assumptions}

We consider the initial boundary value problem governed by the wave equation in a domain $\Omega \subset \mathbb{R}^d$, for $d=2,3$, with smooth boundary $\partial \Omega$. The pressure field satisfies,
\numparts 
\begin{eqnarray} 
\ddot{u} - \C^{2} \Delta u = 0 \quad &&\text{in $(0,T) \times \Omega$}, \label{Eqn.MainIBVP01} \\
\partial_{\nu} u + \gamma \dot{u} = 0 \quad &&\text{on $(0,T) \times \partial \Omega$} \label{Eqn.MainIBVP02}, \\
u = u_{0} \quad \text{and} \quad \dot{u} = u_{1} \quad &&\text{on $\{t=0 \} \times \Omega$} \label{Eqn.MainIBVP03},
\end{eqnarray} 
\endnumparts
on a sufficiently large window of time $(0,T)$ where $0 < T < \infty$. The impedance $\gamma > 0$ models the presence of partially absorbing ultrasound sensors on the boundary $\partial \Omega$ and $\partial_{\nu}$ denotes the outward normal derivative. Consider also the forward mapping given by
\begin{eqnarray} 
\Lambda_{\C} (u_{0},u_{1}) := u|_{(0,T) \times \partial \Omega } \label{Eqn.MeasMap}
\end{eqnarray}
where $u$ solves the system (\ref{Eqn.MainIBVP01})-(\ref{Eqn.MainIBVP03}). In the PAT scenario, it is common to assume that $u_{1} = 0$. However, the mathematical analysis allows us to consider non-vanishing $u_{1}$.

We also consider a possibly different media characterized by a wave speed $\tilde{\C}$ with corresponding wave field $\tilde{u}$ that satisfies,
\numparts 
\begin{eqnarray} 
\ddot{\tilde{u}} - \tilde{\C}^{2} \Delta \tilde{u}  = 0 \quad &\text{in $(0,T) \times \Omega$}, \label{Eqn.tildeIBVP01} \\
\partial_{\nu} \tilde{u} + \gamma \dot{\tilde{u}} = 0 \quad &\text{on $(0,T) \times \partial \Omega$} \label{Eqn.tildeIBVP02}, \\
\tilde{u} = \tilde{u}_{0} \quad \text{and} \quad \dot{\tilde{u}} = \tilde{u}_{1} \quad &\text{on $\{t=0 \} \times \Omega$} \label{Eqn.tildeIBVP03}.
\end{eqnarray} 
\endnumparts
The wave speed $\tilde{\C}$ induces the definition of the corresponding forward map,
\begin{eqnarray} 
\Lambda_{\tilde{\C}} (\tilde{u}_{0},\tilde{u}_{1}) := \tilde{u}|_{(0,T) \times \partial \Omega} \label{Eqn.MeasMapTilde}
\end{eqnarray}
where $\tilde{u}$ solves the system (\ref{Eqn.tildeIBVP01})-(\ref{Eqn.tildeIBVP03}). Notice that the boundary conditions (\ref{Eqn.MainIBVP02}) and (\ref{Eqn.tildeIBVP02}) are the same, which is implied if the media properties on $\partial \Omega$ are identical for both problems.

The goal of this paper is to provide reasonable conditions on the wavespeed $\C$, the initial state $(u_{0},u_{1})$, the domain $\Omega$ and time $T$ to guarantee that the data $\Lambda_{\C}(u_{0},u_{1})$ determines the triplet $(\C,u_{0},u_{1})$ uniquely. This cannot be done in general. Therefore, we  restrict our attention to a problem satisfying certain conditions that we list as follows.

\begin{assumption}[] \label{Assump.Speeds} 
Let the wave speeds $\C$ and $\tilde{\C}$ be smooth in $\clo{\Omega}$ and the impedance $\gamma$ be positive and smooth in $\partial \Omega$. Moreover, let
\begin{eqnarray}
\underline{\C} \leq \C(x)  \quad \text{and} \quad \underline{\C} \leq \tilde{\C}(x)  \qquad \text{for all $x \in \Omega$} 
\label{Eqn.Assumption00speeds}
\end{eqnarray}
for some $ \underline{\C} > 0$. Also assume that
\begin{eqnarray}
 \| \C \|_{W^{1,\infty}(\Omega)} \leq \overline{\C} \quad \text{and} \quad  \| \tilde{\C} \|_{W^{1,\infty}(\Omega)} \leq \overline{\C}
\label{Eqn.Assumption00smooth}
\end{eqnarray}
for some $0 < \overline{\C} < \infty$. Here $W^{1,\infty}(\Omega)$ is the Sobolev space defined by
\begin{eqnarray*}
W^{1,\infty}(\Omega) = \left\{ v \in L^{\infty}(\Omega) ~:~ |\nabla v| \in L^{\infty}(\Omega)  \right\}.
\end{eqnarray*}
\end{assumption}

The next assumption is a geometric condition that ensures the observability of waves from the boundary. This is the so--called \textit{geometric control} or \textit{non--trapping} condition from Bardos--Lebeau--Rauch \cite{Bardos1992}. See also \cite{Triggiani2002,Gulliver2004,Lions1988,Tucsnak2009,GlowinskiLionsHe2008,Las-Tri-2000,Burq1997a,Burq1997b} for references on controllability and observability theory for hyperbolic equations.

\begin{assumption}[Non--trapping condition] \label{Assump.GCC} 
Let $\Omega$ be a simply--connected domain with smooth boundary $\partial \Omega$. For the Riemannian manifold $(\Omega, \C^{-2} dx^2)$, let geodesic rays have finite--order contact with the boundary. Assume there exists $T_{\rm o} < \infty$ such that any geodesic ray originating from any point in $\Omega$ at $t=0$, reaches $\partial \Omega$ at a non--diffractive point before time $t=T_{\rm o}$. Let $T > T_{\rm o}$.
\end{assumption}

It will become clear that we will rely on the following restrictions on the unknown initial state of the pressure field.

\begin{assumption}[] \label{Assump.Rel} 
Let $u_{0} , \tilde{u}_{0} \in H_{0}^{1}(\Omega)$ and $u_{1}, \tilde{u}_{1} \in H^{0}(\Omega)$. Let the following bounds hold
\begin{eqnarray}
\| \nabla u_{0} \|^2_{H^{0}(\Omega)} + \| u_{1} \|^2_{H^{0}(\Omega)} \leq K \quad \text{and} \quad
\| \nabla \tilde{u}_{0} \|^2_{H^{0}(\Omega)} + \| \tilde{u}_{1} \|^2_{H^{0}(\Omega)}  \leq K 
 \label{Eqn.AssumptionBounds02a} \\ 
k \leq \| u_{0} \|^2_{H^{0}(\Omega)} + \| u_{1} \|^2_{H^{-1}(\Omega)}  \quad \text{and} \quad
k \leq \| \tilde{u}_{0} \|^2_{H^{0}(\Omega)} + \| \tilde{u}_{1} \|^2_{H^{-1}(\Omega)} 
\label{Eqn.AssumptionBounds02b}
\end{eqnarray}
for some constants $0 < k < K < \infty$.
\end{assumption}

Before we state the main result, we wish to comment on the relevance of these three assumptions. In the context of PAT, the Assumption \ref{Assump.Speeds} is quite reasonable because, even if the actual wave speed is unknown, lower and upper bounds are readily available due to the nature of biological tissues. In other words, the types of tissue are known (muscular, granular, stromal, cancerous or fatty tissues, and blood or cerebrospinal fluid), but not their distribution within the domain of interest.

Assumption \ref{Assump.GCC} is essential from the physical and mathematical point of view, since it ensures that the energy of the unknown initial pressure reaches the boundary, in finite time, where it can be measured. This assumption allows us to observe, in a stable manner, the initial state of the pressure field from the boundary. 

Assumption \ref{Assump.Rel} 
is verifiable for photoacoustic imaging because $\tilde{u}_{1} = 0$ and the initial pressure profile is given by
\begin{eqnarray}
u_{0}(x) = G(x) \sigma(x) I(x) 
\end{eqnarray}
where $G$ is the Gr\"{u}neisen coefficient, $\sigma$ is the optical absorption coefficient of the medium and $I$ is the intensity of the probing light. The first two factors have well-known upper bounds in biological tissues. The intensity of light is chosen at the boundary and satisfies an elliptic equation in the diffusive regime. Therefore, it is also bounded above. This means that it is possible to find a finite constant $K$ for bound (\ref{Eqn.AssumptionBounds02a}). Now, in order for PAT to work, enough electromagnetic energy must be absorbed by the medium to trigger a measurable acoustic wave. Hence, as long as the Gr\"{u}neisen and absorption coefficients and the intensity of light are bounded away from zero, then it is possible to find a non-vanishing constant $k$ for bound (\ref{Eqn.AssumptionBounds02b}). The discrepancy in the norms employed in bounds (\ref{Eqn.AssumptionBounds02a}) and (\ref{Eqn.AssumptionBounds02b}) is a technicality that we were not able to avoid. In brief, Assumption \ref{Assump.Rel} means that even though the initial pressure state is unknown, its energy has known upper and lower bounds.

Under these conditions, we obtain the main result of this paper.

\begin{theorem}[\textbf{Main Result}] \label{Thm.Main1}
Let Assumptions \ref{Assump.Speeds} and \ref{Assump.GCC} hold for the domain $\Omega$, the wave speeds $\C$ and $\tilde{\C}$, and the time $T_{\rm o} < T < \infty$. Let Assumption \ref{Assump.Rel} hold for the initial states of the pressure fields. There exist positive $\epsilon = \epsilon(\Omega,\C,\underline{\C},\overline{\C},\gamma,k,K)$ so that if
\begin{eqnarray}
\frac{\| \C^{-2} - \tilde{\C}^{-2} \|^{2}_{W^{1,\infty}(\Omega)}}{\| \C^{-2} \|^{2}_{W^{1,\infty}(\Omega)}} & \leq  \epsilon \, \frac{ \| u_{0} -  \tilde{u}_{0} \|^{2}_{H^{0}(\Omega)} + \| u_{1} - \tilde{u}_{1} \|^{2}_{H^{-1}(\Omega)} }{   \| u_{0} \|^{2}_{H^{0}(\Omega)} + \| u_{1} \|^{2}_{H^{-1}(\Omega)}  } \label{Eqn.AssumptionSection}
\end{eqnarray}
then
\begin{eqnarray} \label{eqn.SectionalStability}
\fl \frac{ \| u_{0} - \tilde{u}_{0} \|^{2}_{H^{0}(\Omega)} + \| u_{1} - \tilde{u}_{1} \|^{2}_{H^{-1}(\Omega)} }{ \| u_{0} \|^{2}_{H^{0}(\Omega)} + \| u_{1} \|^{2}_{H^{-1}(\Omega)} }  \leq \frac{C}{T} \frac{ \| \Lambda_{\C}(u_{0},u_{1}) - \Lambda_{\tilde{\C}}(\tilde{u}_{0},\tilde{u}_{1}) \|^{2}_{H^{1}((0,T) ; H^{0}(\partial \Omega))} }{ \| \Lambda_{\C}(u_{0},u_{1}) \|^{2}_{H^{0}((0,T) \times \partial \Omega)}  }
\end{eqnarray}
where $C = C(\Omega,\C,\underline{\C},\overline{\C},\gamma)$. In particular, $\Lambda_{\C}(u_{0},u_{1}) = \Lambda_{\tilde{\C}}(\tilde{u}_{0},\tilde{u}_{1})$ implies that $ u_{0} = \tilde{u}_{0}$ and $u_{1} = \tilde{u}_{1}$.
\end{theorem}

\begin{figure}[!h]
\begin{center}
  \includegraphics[width=0.60 \textwidth]{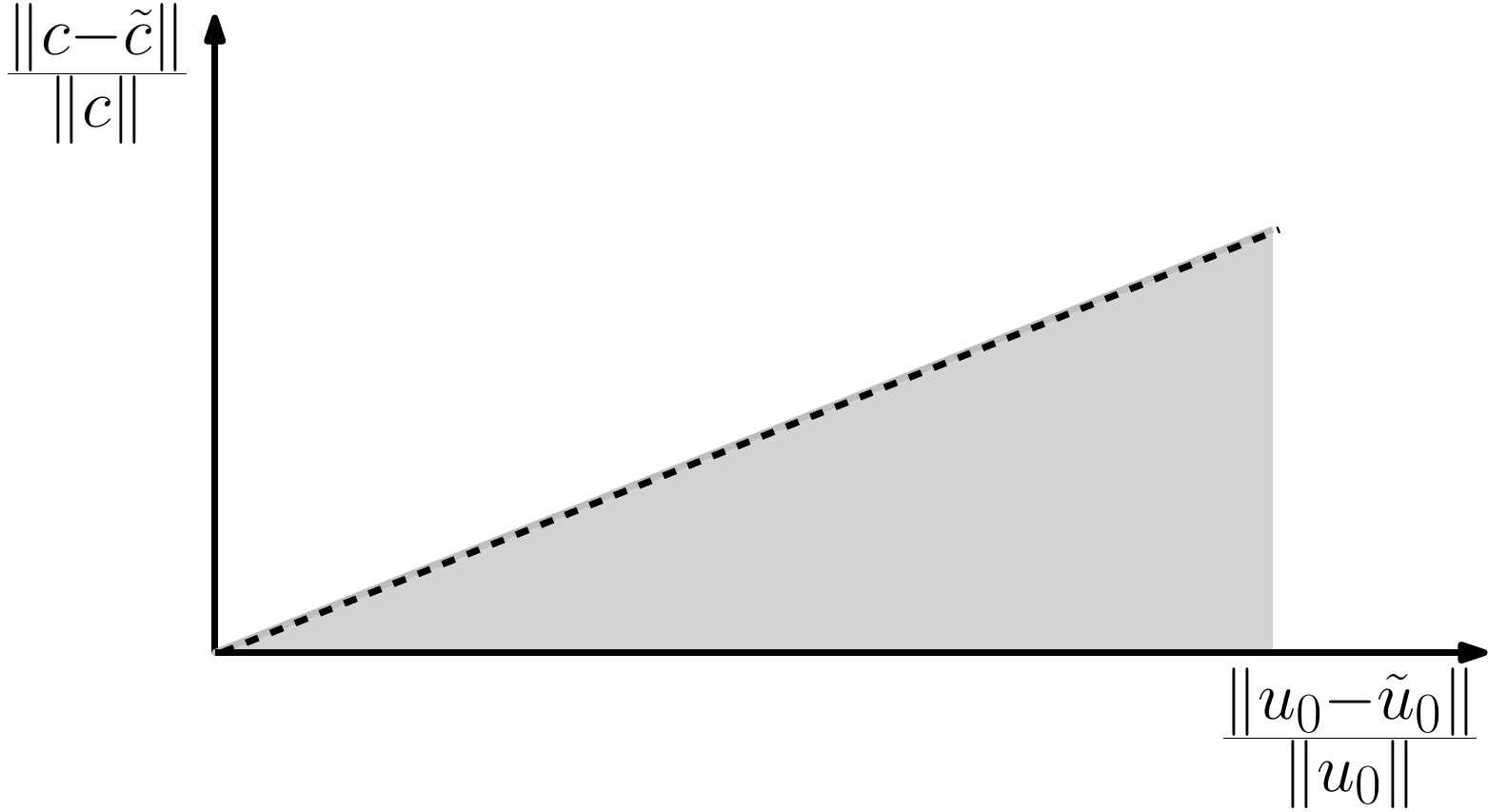}\\
  \caption{Illustration of condition (\ref{Eqn.AssumptionSection}) for the special case $u_{1} = \tilde{u}_{1} = 0$. For a fixed pair $(u_{0},c)$, then a different pair $(\tilde{u}_{0}, \tilde{c})$ from the pre--image of the shaded region cannot produce the same boundary measurements. The norms are understood as in Theorem \ref{Thm.Main1}. The slope of the tilted line defining the shaded region is $\epsilon^{1/2} > 0$.} \label{Figure.Sketch}
  \end{center}
\end{figure}

In the context of PAT, where $u_{1} = \tilde{u}_{1} = 0$, this theorem states that given two different pairs $(\C,u_{0})$ and $(\tilde{\C},\tilde{u}_{1})$, if the relative difference between the wave speeds is much smaller than the relative difference between the initial pressure profiles, then these two pairs of data cannot induce the same acoustic measurements at the boundary. Condition (\ref{Eqn.AssumptionSection}) is illustrated in Figure \ref{Figure.Sketch}. 

Notice that Theorem \ref{Thm.Main1} is not a statement of uniqueness in the usual sense. The estimate (\ref{eqn.SectionalStability}) is only valid under condition (\ref{Eqn.AssumptionSection}). This assumption does not include a neighborhood of the triplet $(\C,u_{0},u_{1})$. It only considers a small \textit{conical region} as illustrated in Figure \ref{Figure.Sketch}. No other triple $(\tilde{\C},\tilde{u}_{0},\tilde{u}_{1})$ in that region can produce the same boundary measurements as the triple $(\C,u_{0},u_{1})$.

\section{Proof of the Main Result} \label{Sec:Proof}
We start by setting up the initial value problem for the contrast $w = u - \tilde{u}$. Notice that $w$ satisfies,
\numparts  
\begin{eqnarray} 
\C^{-2} \ddot{w} - \Delta w = f    \qquad &\text{in $(0,T) \times \Omega$}, \label{Eqn.diff_wave01} \\
 w = m \qquad &\text{on $(0,T) \times \partial \Omega$} \label{Eqn.diff_wave02}, \\
w = w_{0} \quad \text{and} \quad \dot{w} = w_{1} \qquad &\text{on $\{t=0 \} \times \Omega$} \label{Eqn.diff_wave03},
\end{eqnarray} 
\endnumparts 
where $f = \left( \tilde{\C}^{-2} - \C^{-2} \right) \ddot{\tilde{u}}$, $w_{0} = u_{0} - \tilde{u}_{0}$, $w_{1} = u_{1} - \tilde{u}_{1}$, and $m = \Lambda_{\C}(u_{0},u_{1}) - \Lambda_{\tilde{\C}}(\tilde{u}_{0},\tilde{u}_{1})$. In terms of regularity, we have that $\ddot{\tilde{u}} \in H^{0}((0,T);H^{-1}(\Omega))$ because $\tilde{u}$ is a weak solution to the wave equation as implied by the regularity of the initial conditions in Assumption \ref{Assump.Rel} \cite{EvansPDE}. Similarly, $m \in H^{1}((0,T)\times \partial \Omega)$ and $m|_{t=0} = 0$ because $\left( u_{0} - \tilde{u}_{0} \right) \in H_{0}^{1}(\Omega)$ as required by. 
Assumption \ref{Assump.Rel}.

By virtue of linearity, we can decompose $w$ as follows $w = w^{(1)} + w^{(2)}$ where 
\numparts 
\begin{eqnarray} 
\C^{-2} \ddot{w}^{(1)} - \Delta w^{(1)} = 0    \qquad &\text{in $(0,T) \times \Omega$}, \label{Eqn.diff_1_wave01} \\
 w^{(1)} = 0 \qquad &\text{on $(0,T) \times \partial \Omega$} \label{Eqn.diff_1_wave02}, \\
w^{(1)} = w_{0} \quad \text{and} \quad \dot{w}^{(1)} = w_{1} \qquad &\text{on $\{t=0 \} \times \Omega$}, \label{Eqn.diff_1_wave03}
\end{eqnarray} 
\endnumparts
and
\numparts 
\begin{eqnarray} 
\C^{-2} \ddot{w}^{(2)} - \Delta w^{(2)} = f    \qquad &\text{in $(0,T) \times \Omega$}, \label{Eqn.diff_2_wave01} \\
 w^{(2)} = m \qquad &\text{on $(0,T) \times \partial \Omega$} \label{Eqn.diff_2_wave02}, \\
w^{(2)} = 0 \quad \text{and} \quad \dot{w}^{(2)} = 0 \qquad &\text{on $\{t=0 \} \times \Omega$}. \label{Eqn.diff_2_wave03}
\end{eqnarray} 
\endnumparts

Now we proceed to state some lemmas concerning these initial boundary value problems. In what follows, the constant $C > 0$ will be a generic constant that change its value from inequality to inequality. However, $C$ does not depend on the profile of $\tilde{\C}$ (only on its lower and upper bounds) and does not depend on the solution to any of the boundary value problems.

The first lemma is a well--known result concerning boundary observability for hyperbolic equations. See classical references \cite{Bardos1992,Lions1988,Burq1997a,Burq1997b} or more recent related works \cite{Tucsnak2009,GlowinskiLionsHe2008,Las-Tri-2000,Triggiani2002,Gulliver2004,Alabau-Boussouira2013a,Alabau-Boussouira2014}.

\begin{lemma}[Boundary Observability] \label{Lemma.001}
Let $w^{(1)}$ solve the initial boundary value problem (\ref{Eqn.diff_1_wave01})-(\ref{Eqn.diff_1_wave03}) and Assumption \ref{Assump.GCC} hold. Then $w^{(1)}$ satisfies the following boundary observability estimate,
\begin{eqnarray*} 
T \left( \| w_{0} \|^{2}_{H^{0}(\Omega)} + \| w_{1} \|^{2}_{H^{-1}(\Omega)} \right) \leq C \| \partial_{\nu} w^{(1)} \|^{2}_{H^{-1}((0,T) \times \partial \Omega)}
\end{eqnarray*} 
for some positive constant $C=C(\Omega,\C)$.
\end{lemma}

The next lemmas are stability results in less-regular spaces for initial boundary value problems for the wave equation using the transposition method \cite{Lasiecka1986}.

\begin{lemma}[] \label{Lemma.002}
Let $w^{(2)}$ solve the initial boundary value problem (\ref{Eqn.diff_2_wave01})-(\ref{Eqn.diff_2_wave03}). Then the following stability estimate holds
\begin{eqnarray*} 
\| \partial_{\nu} w^{(2)} \|^{2}_{H^{-1}((0,T) \times \partial \Omega)} \leq C \left( \| m \|^{2}_{H^{0}((0,T) \times \partial \Omega)} + \| f \|^{2}_{H^{0}((0,T) ; H^{-1}(\Omega))} \right)
\end{eqnarray*} 
for some constant $C = C(\Omega,\C)$.
\end{lemma}

\begin{lemma}[] \label{Lemma.003}
Let $u$ solve the initial boundary value problem (\ref{Eqn.MainIBVP01})-(\ref{Eqn.MainIBVP03}) and $\Lambda_{\C}(u_{0},u_{1})$ be given by (\ref{Eqn.MeasMap}). Then the following stability estimate holds
\begin{eqnarray*} 
\| \Lambda_{\C}(u_{0},u_{1}) \|^{2}_{H^{0}((0,T) \times \partial \Omega)} \leq C \left(  \| u_{0} \|^{2}_{H^{0}(\Omega) }  +  \| u_{1} \|^{2}_{H^{-1}(\Omega)} \right)
\end{eqnarray*} 
for some constant $C = C(\Omega,\C,\gamma)$.
\end{lemma}

\label{Eqn.MainIBVP}

With the above lemmas, we are ready to prove the main result of this paper.

\begin{proof}[\textbf{Proof of Theorem \ref{Thm.Main1}}]
In order to apply Lemmas \ref{Lemma.001}--\ref{Lemma.002}, it only remains to estimate the norm of $f$. Recall that $f = \left( \tilde{\C}^{-2} - \C^{-2} \right) \ddot{\tilde{u}}$. Let 
\begin{eqnarray*}
g = \frac{\tilde{\C}^{-2} - \C^{-2}}{ \tilde{\C}^{-2} }
\end{eqnarray*}
so that $f = g \tilde{\C}^{-2} \ddot{\tilde{u}}$. Now take $v \in H^{0}((0,T);H_{0}^{1}(\Omega))$ with $\| v \|_{H^{0}((0,T);H_{0}^{1}(\Omega))} \leq 1$ and consider that
\begin{eqnarray*}
 | \la f, v \ra| &= | \la \tilde{\C}^{-2} \ddot{\tilde{u}}, g v \ra| \leq |\la \nabla \tilde{u} , \nabla (g v) \ra_{H^{0}((0,T) ; H^{0}(\Omega))}| \\
&\leq \| \nabla \tilde{u} \|_{H^{0}((0,T);H^{0}(\Omega))}  \left(  \| v \nabla g \|_{H^{0}((0,T);H^{0}(\Omega))} + \| g \nabla v \|_{H^{0}((0,T);H^{0}(\Omega))} \right) \\
&\leq \| g \|_{W^{1,\infty}(\Omega)} \| \nabla \tilde{u} \|_{H^{0}((0,T);H^{0}(\Omega))} \| v \|_{H^{0}((0,T);H_{0}^{1}(\Omega))} \\
&\leq \| g \|_{W^{1,\infty}(\Omega)} \| \nabla \tilde{u} \|_{H^{0}((0,T);H^{0}(\Omega))}
\end{eqnarray*}
for $v$ is arbitrary. Using the standard energy estimate \cite{EvansPDE,Lasiecka1986} for the term $\nabla \tilde{u}$ and the assumed bound (\ref{Eqn.AssumptionBounds02a}), we obtain
\begin{eqnarray}
\fl \| f \|^{2}_{H^{0}((0,T);H^{-1}(\Omega))} \leq C T \| g \|^{2}_{W^{1,\infty}(\Omega)}   \left( \| \nabla \tilde{u}_{0} \|^{2}_{H^{0}(\Omega)} + \| \tilde{u}_{1} \|^{2}_{H^{0}(\Omega)}  \right)  \leq K  CT\| g \|^{2}_{W^{1,\infty}(\Omega)} \label{Eqn.Bound_f}
\end{eqnarray}
where $C=C(\Omega,\underline{\C},\overline{\C})$. Now, after some algebra, we obtain that
\begin{eqnarray}
\| g \|_{W^{1,\infty}(\Omega)} \leq \left( 1 +  \| \tilde{\C} \|^{2}_{L^{\infty}(\Omega)} \| \tilde{\C}^{-2}  \|_{W^{1,\infty}(\Omega)}   \right) \| \tilde{\C} \|^{2}_{L^{\infty}(\Omega)} \| \tilde{\C}^{-2} - \C^{-2} \|_{W^{1,\infty}(\Omega)} .  \nonumber
\end{eqnarray}
From the assumed bounds (\ref{Eqn.Assumption00speeds}) and (\ref{Eqn.Assumption00smooth}), we obtain a constant $C=C(\underline{\C},\overline{\C})$ such that $\| \C^{-2} \|_{W^{1,\infty}(\Omega)} \leq C$ and $\| \tilde{\C}^{-2} \|_{W^{1,\infty}(\Omega)} \leq C$. Consequently, there is another constant $C=C(\underline{\C},\overline{\C})$ such that,
\begin{eqnarray}
\| g \|_{W^{1,\infty}(\Omega)} \leq C \frac{\| \C^{-2} - \tilde{\C}^{-2} \|_{W^{1,\infty}(\Omega)}}{\| \C^{-2} \|_{W^{1,\infty}(\Omega)}}  \label{Eqn.Bound_g}
\end{eqnarray}
Plugging (\ref{Eqn.Bound_g}) into (\ref{Eqn.Bound_f}) and using the assumed bound (\ref{Eqn.AssumptionSection}), we obtain
\begin{eqnarray}
\| f \|^{2}_{H^{0}((0,T) ; H^{-1}(\Omega))} \leq  \epsilon K C T  \frac{  \| w_{0} \|^{2}_{H^{0}(\Omega) }  +  \| w_{1} \|^{2}_{H^{-1}(\Omega)} }{ \| u_{0} \|^{2}_{H^{0}(\Omega) }  +  \| u_{1} \|^{2}_{H^{-1}(\Omega)} }  \label{Eqn.Bound003}
\end{eqnarray}
where $C=C(\Omega,\underline{\C},\overline{\C})$.

Now, combining Lemmas \ref{Lemma.001} and \ref{Lemma.002}, the estimate (\ref{Eqn.Bound003}), and the decomposition $w = w^{(1)} + w^{(2)}$, we obtain
\begin{eqnarray*}
\fl \| w_{0} \|^{2}_{H^{0}(\Omega)} + \| w_{1} \|^{2}_{H^{-1}(\Omega)} & \leq \frac{C}{T} \| \partial_{\nu} w^{(1)} \|^{2}_{H^{-1}((0,T) \times \partial \Omega)} \\ & \leq \frac{C}{T} \left( \| \partial_{\nu} w \|^{2}_{H^{-1}((0,T) \times \partial \Omega)} + \| \partial_{\nu} w^{(2)} \|^{2}_{H^{-1}((0,T) \times \partial \Omega)}  \right) \\
& \leq \frac{C}{T} \left( \| \gamma \dot{m} \|^{2}_{H^{0}((0,T) \times \partial \Omega)}  + \| m \|^{2}_{H^{0}((0,T) \times \partial \Omega)} + \| f \|^{2}_{H^{0}((0,T) ; H^{-1}(\Omega))}  \right) \\
& \leq \frac{C}{T} \left(  \| m \|^{2}_{H^{1}((0,T) ; H^{0}(\partial \Omega))} + \epsilon K T \frac{  \| w_{0} \|^{2}_{H^{0}(\Omega) }  +  \| w_{1} \|^{2}_{H^{-1}(\Omega)} }{ \| u_{0} \|^{2}_{H^{0}(\Omega) }  +  \| u_{1} \|^{2}_{H^{-1}(\Omega)} }  \right)
\end{eqnarray*}
where $C=C(\Omega,\C,\underline{\C},\overline{\C},\gamma)$, and $\partial_{\nu} w = - \gamma \dot{w}$ because $w = u - \tilde{u}$ where $u$ and $\tilde{u}$ satisfy (\ref{Eqn.MainIBVP02}) and (\ref{Eqn.tildeIBVP02}), respectively, on the boundary. Rearranging some terms and using Lemma \ref{Lemma.003} and the assumed bounds (\ref{Eqn.AssumptionBounds02a})-(\ref{Eqn.AssumptionBounds02b}), we get
\begin{eqnarray*}
 \left( \frac{k}{K} - \epsilon C   \right) \frac{  \| w_{0} \|^{2}_{H^{0}(\Omega) }  +  \| w_{1} \|^{2}_{H^{-1}(\Omega)} }{  \| u_{0} \|^{2}_{H^{0}(\Omega) }  +  \| u_{1} \|^{2}_{H^{-1}(\Omega)}  } 
 \leq \frac{C}{T} \frac{ \| m \|^{2}_{H^{1}((0,T); H^{0}(\partial \Omega))} }{ \| \Lambda_{\C}(u_{0},u_{1}) \|^{2}_{H^{0}((0,T) \times \partial \Omega)}  }
\end{eqnarray*}
where $C=C(\Omega,c,\underline{\C},\overline{\C},\gamma)$ does not depend on $\epsilon$ or $\tilde{\C}$. Therefore, for a choice $\epsilon < k /\left( K C \right)$, we obtain the desired result.
\end{proof}

\section{Final Remarks} \label{Sec:Final}

We end by offering some remarks concerning the practical significance of the theoretical results from the previous sections. In particular, we are interested in the successful design of iterative algorithms for the joint reconstruction of the pressure profile and wave speed. The main concern is to ensure that iterates remain within the region of uniqueness defined by (\ref{Eqn.AssumptionSection}).

We seek to recover $(u_{0},\C)$ and assume that $u_{1} = 0$ as in the case in PAT. Let $(u^{(n)},\C^{(n)})$ for $n=0,1,2, ... $ be a sequence of iterates converging to $(u_{0},\C)$. Many of the iterative algorithms display \emph{linear} convergence (such as Banach fixed-point, gradient descent, Landweber or conjugate gradient iterations \cite{Atkinson-Han-Book-2001,Stoer2002}), meaning that the error behaves like
\numparts \label{Eqn.linear_conv}
\begin{eqnarray}
\alpha_{u} \| u_{0} - u^{(n)} \| &\leq \| u_{0} - u^{(n+1)} \| \leq \beta_{u} \| u_{0} - u^{(n)} \|    \label{Eqn.linear_conv_u} \\
\alpha_{\C} \| \C - \C^{(n)} \| &\leq \| \C - \C^{(n+1)} \| \leq \beta_{\C} \| \C - \C^{(n)} \| \label{Eqn.linear_conv_c}   
\end{eqnarray}
\endnumparts
for $n=0,1,2, ...$ and for some constants $0 < \alpha_{u} < \beta_{u} < 1$ and $0 < \alpha_{\C} < \beta_{\C} < 1$, in the appropriate norms. Using the estimates (\ref{Eqn.linear_conv_u})-(\ref{Eqn.linear_conv_c}), then we arrive at
\begin{eqnarray}
\frac{\| \C - \C^{(n+1)} \|}{ \| \C \| } \leq \epsilon^{1/2} \left( \frac{\beta_{\C}}{\alpha_{u}} \right) \frac{\| u_{0} - u^{(n+1)} \|}{ \| u_{0} \| } ,
\end{eqnarray}
provided that $(u^{(n)},\C^{(n)})$ satisfies (\ref{Eqn.AssumptionSection}). Hence, the next iterate $(u^{(n+1)},\C^{(n+1)})$ remains in the region of uniqueness defined by (\ref{Eqn.AssumptionSection}) if
\begin{eqnarray}
\beta_{\C} \leq \alpha_{u}.  \label{Eqn.alpha_beta}
\end{eqnarray}
We can show inductively that all the iterates $(u^{(n)},\C^{(n)})$ for $n=1,2,3 ,...$ satisfy    (\ref{Eqn.AssumptionSection}) provided that (\ref{Eqn.alpha_beta}) holds and that the initial guess $(u^{(0)},\C^{(0)})$ satisfies (\ref{Eqn.AssumptionSection}).

In practical terms this means that for a given linear iterative algorithm, the iterates $u^{(n)}$ for the pressure profile must be relaxed so that they do not converge faster than the iterates $\C^{(n)}$ for the wave speed. Or alternatively, the convergence for the wave speed iterates should be accelerated. Acceleration of linear convergence can be achieved by methods of Nesterov \cite{Nesterov2018} or Aitken \cite[\S 5.10]{Stoer2002}. Acceleration can also be obtained with additional constraints on the wave speed. For instance, if the wave speed is assumed to belong to finite-dimensional parametric spaces, then the compactness of the projection accelerates methods such as the conjugate gradient or Landweber methods. In the computational setting, both $u_{0}$ and $\C$ belong to finite-dimensional spaces following the discretization of the domain $\Omega$ and governing differential equations. In that case, the parametrization of $\C$ should belong to a space with much lower dimension than that of the parametrization of $u_{0}$. This can be achieved with a two-mesh approach, one mesh being much coarser than the other. Another approach was employed in \cite{Matthews2018} where each pixel value of the wave speed was assumed to belong to one of a few tissue types with each tissue type having a uniform sound speed.

\section*{Acknowledgements} \label{Sec:Ack}
This work was partially supported by NSF grant DMS-1712725. The author would like to thank Texas Children's Hospital for its support and for the research-oriented environment provided by the Predictive Analytics Laboratory.


\section*{References}

\bibliographystyle{siamplain}
\bibliography{InvProblemBiblio}

\begin{thebibliography}{10}

\bibitem{Acosta-Montalto-2015}
{\sc S.~Acosta and C.~Montalto}, {\em {Multiwave imaging in an enclosure with
  variable wave speed}}, Inverse Probl., 31 (2015), p.~065009,
  \url{https://doi.org/10.1088/0266-5611/31/6/065009}.

\bibitem{Acosta-Montalto-2016}
{\sc S.~Acosta and C.~Montalto}, {\em {Photoacoustic imaging taking into
  account thermodynamic attenuation}}, Inverse Probl., 32 (2016), p.~115001,
  \url{https://doi.org/10.1088/0266-5611/32/11/115001}.

\bibitem{Alabau-Boussouira2014}
{\sc F.~Alabau-Boussouira}, {\em {Insensitizing exact controls for the scalar
  wave equation and exact controllability of 2-coupled cascade systems of PDE's
  by a single control}}, Math. Control. Signals, Syst., 26 (2014), pp.~1--46,
  \url{https://doi.org/10.1007/s00498-013-0112-8}.

\bibitem{Alabau-Boussouira2013a}
{\sc F.~Alabau-Boussouira and M.~Leautaud}, {\em {Indirect controllability of
  locally coupled wave-type systems and applications}}, J. des Math. Pures
  Appl., 99 (2013), pp.~544--576,
  \url{https://doi.org/10.1016/j.matpur.2012.09.012}.

\bibitem{Atkinson-Han-Book-2001}
{\sc K.~Atkinson and W.~Han}, {\em {Theoretical Numerical Analysis: A
  Functional Analysis Framework}}, vol.~39 of Texts in Applied Mathematics, New
  York : Springer-Verlag, 2001.

\bibitem{Bardos1992}
{\sc C.~Bardos, G.~Lebeau, and J.~Rauch}, {\em {Sharp Sufficient Conditions for
  the Observation, Control, and Stabilization of Waves from the Boundary}},
  SIAM J. Control Optim., 30 (1992), pp.~1024--1065,
  \url{https://doi.org/10.1137/0330055}.

\bibitem{Beard2011}
{\sc P.~Beard}, {\em {Biomedical photoacoustic imaging}}, Interface Focus, 1
  (2011), pp.~602--31, \url{https://doi.org/10.1098/rsfs.2011.0028}.

\bibitem{Belhachmi2016}
{\sc Z.~Belhachmi, T.~Glatz, and O.~Scherzer}, {\em {A direct method for
  photoacoustic tomography with inhomogeneous sound speed}}, Inverse Probl., 32
  (2016), p.~045005, \url{https://doi.org/10.1088/0266-5611/32/4/045005}.

\bibitem{Burq1997b}
{\sc N.~Burq}, {\em {Contr{\^{o}}labilit{\'{e}} exacte des ondes dans des
  ouverts peu r{\'{e}}guliers}}, Asymptot. Anal., 14 (1997), pp.~157--191,
  \url{https://doi.org/10.3233/ASY-1997-14203}.

\bibitem{Burq1997a}
{\sc N.~Burq and P.~Gerard}, {\em {A necessary and sufficient condition for the
  exact controllability of the wave equation}}, C. R. Acad. Sci. Paris, 325
  (1997), pp.~749--752.

\bibitem{Chervova2016}
{\sc O.~Chervova and L.~Oksanen}, {\em {Time reversal method with stabilizing
  boundary conditions for photoacoustic tomography}}, Inverse Probl., 32
  (2016), p.~125004, \url{https://doi.org/10.1088/0266-5611/32/12/125004}.

\bibitem{Cox2012}
{\sc B.~Cox, J.~G. Laufer, S.~R. Arridge, and P.~C. Beard}, {\em {Quantitative
  spectroscopic photoacoustic imaging: A review}}, J. Biomed. Opt., 17 (2012),
  pp.~061202--1--061202--22, \url{https://doi.org/10.1117/1.JBO.17.6.061202}.

\bibitem{Dean-Ben2014}
{\sc X.~L. Dean-Ben, V.~Ntziachristos, and D.~Razansky}, {\em {Effects of small
  variations of speed of sound in optoacoustic tomographic imaging}}, Med.
  Phys., 41 (2014), p.~073301, \url{https://doi.org/10.1118/1.4875691}.

\bibitem{EvansPDE}
{\sc L.~C. Evans}, {\em {Partial Differential Equations}}, vol.~19 of Graduate
  Studies in Mathematics, American Mathematical Society, Providence, RI, 1998.

\bibitem{Finch2004}
{\sc D.~Finch, S.~Patch, and Rakesh}, {\em {Determining a function from its
  mean values over a family of spheres}}, SIAM J. Math. Anal., 35 (2004),
  pp.~1213--1240, \url{https://doi.org/10.1137/S0036141002417814}.

\bibitem{GlowinskiLionsHe2008}
{\sc R.~Glowinski, J.-L. Lions, and J.~He}, {\em {Exact and Approximate
  Controllability for Distributed Parameter Systems : A Numerical Approach}},
  vol.~117 of Encyclopedia of Mathematics and its Applications, Cambridge
  University Press, 2008, \url{https://doi.org/10.1017/CBO9780511721595}.

\bibitem{Gulliver2004}
{\sc R.~Gulliver, I.~Lasiecka, W.~Littman, and R.~Triggiani}, {\em {The case
  for differential geometry in the control of single and coupled PDEs: The
  structural acoustic chamber}}, in Geom. Methods Inverse Probl. PDE Control,
  C.~Croke, I.~Lasiecka, G.~Uhlmann, and M.~Vogelius, eds., vol.~137 of IMA
  volumes in Mathematics and its Applications, Springer-Verlag, New York, 2004,
  pp.~73--181, \url{https://doi.org/10.1007/978-1-4684-9375-7_5}.

\bibitem{Haltmeier2014}
{\sc M.~Haltmeier}, {\em {Universal inversion formulas for recovering a
  function from spherical means}}, SIAM J. Math. Anal., 46 (2014),
  pp.~214--232, \url{https://doi.org/10.1137/120881270}.

\bibitem{Haltmeier2017b}
{\sc M.~Haltmeier, R.~Kowar, and L.~V. Nguyen}, {\em {Iterative methods for
  photoacoustic tomography in attenuating acoustic media}}, Inverse Probl., 33
  (2017), p.~115009, \url{https://doi.org/10.1088/1361-6420/aa8cba}.

\bibitem{Huang2013}
{\sc C.~Huang, K.~Wang, L.~Nie, L.~Wang, and M.~Anastasio}, {\em {Full-wave
  iterative image reconstruction in photoacoustic tomography with acoustically
  inhomogeneous media}}, IEEE Trans. Med. Imaging, 32 (2013), pp.~1097--1110,
  \url{https://doi.org/10.1109/TMI.2013.2254496}.

\bibitem{Huang2016a}
{\sc C.~Huang, K.~Wang, R.~W. Schoonover, L.~V. Wang, M.~A. Anastasio,
  C.~Huang, L.~V. Wang, and K.~Wang}, {\em {Joint reconstruction of absorbed
  optical energy density and sound speed distributions in photoacoustic
  computed tomography: A numerical investigation}}, IEEE Trans. Comput.
  Imaging, 2 (2016), pp.~136--149,
  \url{https://doi.org/10.1109/TCI.2016.2523427}.

\bibitem{Jin2006}
{\sc X.~Jin and L.~V. Wang}, {\em {Thermoacoustic tomography with correction
  for acoustic speed variations}}, Phys. Med. Biol., 51 (2006), pp.~6437--6448,
  \url{https://doi.org/10.1088/0031-9155/51/24/010}.

\bibitem{Kirsch2012}
{\sc A.~Kirsch and O.~Scherzer}, {\em {Simultaneous reconstructions of
  absorption density and wave speed with photoacoustic measurements}}, SIAM J.
  Appl. Math., 72 (2012), pp.~1508--1523,
  \url{https://doi.org/10.1137/110849055}.

\bibitem{Kuchment2008}
{\sc P.~Kuchment and L.~Kunyansky}, {\em {Mathematics of thermoacoustic
  tomography}}, Eur. J. Appl. Math., 19 (2008), pp.~191--224,
  \url{https://doi.org/10.1017/S0956792508007353}.

\bibitem{Kunyansky2007}
{\sc L.~Kunyansky}, {\em {Explicit inversion formulas for the spherical mean
  Radon transform}}, Inverse Probl., 23 (2007), pp.~373--383,
  \url{https://doi.org/10.1088/0266-5611/23/1/021}.

\bibitem{Kunyansky2011}
{\sc L.~Kunyansky}, {\em {Reconstruction of a function from its spherical
  (circular) means with the centers lying on the surface of certain polygons
  and polyhedra}}, Inverse Probl., 27 (2011), p.~025012,
  \url{https://doi.org/10.1088/0266-5611/27/2/025012}.

\bibitem{Lasiecka1986}
{\sc I.~Lasiecka, J.-L. Lions, and R.~Triggiani}, {\em {Nonhomogeneous boudary
  value problems for second order hyperbolic equations}}, J. Math. Pures Appl.,
  65 (1986), pp.~149--192.

\bibitem{Las-Tri-2000}
{\sc I.~Lasiecka and R.~Triggiani}, {\em {Control theory for partial
  differential equations : Continuous and approximation theories}}, vol.~74-75
  of Encyclopedia of mathematics and its applications, Cambridge University
  Press, Cambridge; New York, 2000,
  \url{https://doi.org/10.1017/CBO9781107340848}.

\bibitem{Lions1988}
{\sc J.-L. Lions}, {\em {Exact controllability, stabilization and perturbations
  for distributed systems}}, SIAM Rev., 30 (1988), pp.~1----68,
  \url{https://doi.org/10.1137/1030001}.

\bibitem{Liu2015}
{\sc H.~Liu and G.~Uhlmann}, {\em {Determining both sound speed and internal
  source in thermo- and photo-acoustic tomography}}, Inverse Probl., 31 (2015),
  pp.~1--11, \url{https://doi.org/10.1088/0266-5611/31/10/105005}.

\bibitem{Matthews2017b}
{\sc T.~P. Matthews and M.~A. Anastasio}, {\em {Joint reconstruction of the
  initial pressure and speed of sound distributions from combined photoacoustic
  and ultrasound tomography measurements}}, Inverse Probl., 33 (2017),
  p.~124002, \url{https://doi.org/10.1088/1361-6420/aa9384}.

\bibitem{Matthews2018}
{\sc T.~P. Matthews, J.~Poudel, L.~Li, L.~V. Wang, and M.~Anastasio}, {\em
  {Parameterized joint reconstruction of the initial pressure and sound speed
  distributions for photoacoustic computed tomography}}, SIAM J. Imaging Sci.,
  11 (2018), pp.~1560--1588, \url{https://doi.org/10.1137/17M1153649}.

\bibitem{Natterer2012}
{\sc F.~Natterer}, {\em {Photo-acoustic inversion in convex domains}}, Inverse
  Probl. Imaging, 6 (2012), pp.~315--320,
  \url{https://doi.org/10.3934/ipi.2012.6.315}.

\bibitem{Nesterov2018}
{\sc Y.~Nesterov}, {\em {Lectures on Convex Optimization}}, vol.~137 of
  Springer Optimization and Its Applications, Springer, 2nd~ed., 2018,
  \url{https://doi.org/10.1007/978-3-319-91578-4}.

\bibitem{Oksanen2014}
{\sc L.~Oksanen and G.~Uhlmann}, {\em {Photoacoustic and thermoacoustic
  tomography with an uncertain wave speed}}, Math. Res. Lett., 21 (2014),
  pp.~1199--1214, \url{https://doi.org/10.4310/MRL.2014.v21.n5.a13},
  \url{https://arxiv.org/abs/1307.1618}.

\bibitem{Palacios2016b}
{\sc B.~Palacios}, {\em {Reconstruction for multi-wave imaging in attenuating
  media with large damping coefficient}}, Inverse Probl., 32 (2016), p.~125008,
  \url{https://doi.org/10.1088/0266-5611/32/12/125008}.

\bibitem{Qian2011}
{\sc J.~Qian, P.~Stefanov, G.~Uhlmann, and H.~Zhao}, {\em {An efficient Neumann
  series-based algorithm for thermoacoustic and photoacoustic tomography with
  variable sound speed}}, SIAM J. Imaging Sci., 4 (2011), pp.~850--883,
  \url{https://doi.org/10.1137/100817280}.

\bibitem{Stefanov2009}
{\sc P.~Stefanov and G.~Uhlmann}, {\em {Thermoacoustic tomography with variable
  sound speed}}, Inverse Probl., 25 (2009), p.~075011,
  \url{https://doi.org/10.1088/0266-5611/25/7/075011}.

\bibitem{Stefanov2013b}
{\sc P.~Stefanov and G.~Uhlmann}, {\em {Instability of the linearized problem
  in multiwave tomography of recovery both the source and the speed}}, Inverse
  Probl. Imaging, 7 (2013), pp.~1367--1377,
  \url{https://doi.org/10.3934/ipi.2013.7.1367},
  \url{https://arxiv.org/abs/1211.6217}.

\bibitem{Stefanov2013}
{\sc P.~Stefanov and G.~Uhlmann}, {\em {Recovery of a source term or a speed
  with one measurement and applications}}, Trans. Am. Math. Soc., 365 (2013),
  pp.~5737--5758.

\bibitem{Stefanov2017a}
{\sc P.~Stefanov and Y.~Yang}, {\em {Multiwave tomography with reflectors:
  Landweber's iteration}}, Inverse Probl. Imaging, 11 (2017), pp.~373--401,
  \url{https://doi.org/10.3934/ipi.2017018}.

\bibitem{Stoer2002}
{\sc J.~Stoer and R.~Bulirsch}, {\em {Introduction to Numerical Analysis}},
  vol.~12 of Texts in Applied Mathematics, Springer, 3rd~ed., 2002,
  \url{https://doi.org/10.1007/978-0-387-21738-3}.

\bibitem{Treeby2011}
{\sc B.~Treeby, T.~Varslot, E.~Zhang, J.~Laufer, and P.~Beard}, {\em {Automatic
  sound speed selection in photoacoustic image reconstruction using an
  autofocus approach}}, J. Biomed. Opt., 16 (2011), p.~090501,
  \url{https://doi.org/10.1117/1.3619139},
  \url{http://www.ncbi.nlm.nih.gov/pubmed/21950905}.

\bibitem{Triggiani2002}
{\sc R.~Triggiani and P.~Yao}, {\em {Carleman estimates with no lower-order
  terms for general Riemann wave equations. Global uniqueness and observability
  in one shot}}, Appl. Math. Optim., 46 (2002), pp.~331--375,
  \url{https://doi.org/10.1007/s00245-002-0751-5}.

\bibitem{Tucsnak2009}
{\sc M.~Tucsnak and G.~Weiss}, {\em {Observation and control for operator
  semigroups}}, Birkh{\"{a}}user Advanced Texts, Springer Science
  $\backslash${\&} Business Media, 2009,
  \url{https://doi.org/10.1007/978-3-7643-8994-9}.

\bibitem{Wang-Anastasio-2011}
{\sc K.~Wang and M.~Anastasio}, {\em {Photoacoustic and Thermoacoustic
  Tomography: Image Formation Principles}}, in Handb. Math. Methods Imaging,
  O.~Scherzer, ed., Springer New York, 2011, pp.~781--815.

\bibitem{Wang-2009}
{\sc L.~V. Wang}, {\em {Photoacoustic Imaging and Spectroscopy}}, Boca Raton :
  CRC, 2009.

\bibitem{Wang2012}
{\sc L.~V. Wang and S.~Hu}, {\em {Photoacoustic tomography: In vivo imaging
  from organelles to organs}}, Science (80-. )., 335 (2012), pp.~1458--1462,
  \url{https://doi.org/10.1126/science.1216210}.

\bibitem{Wang-Wu-2007}
{\sc L.~V. Wang and H.~Wu}, {\em {Biomedical optics: Principles and Imaging}},
  Hoboken, N.J. : Wiley-Interscience, 2007,
  \url{https://doi.org/10.1002/9780470177013}.

\end{thebibliography}

\end{document}